\documentclass[12pt]{article}
\def\Box{{\setlength{\unitlength}{1.1 ex}
\begin{picture}(1,1)(-0.2,-0.2)
\put(0,0){\framebox(1,1){}}
\end{picture}}}

\font\elevenbb=msbm10 at 10.95pt

\def\N{\hbox{\elevenbb N}}
\def\M{\hbox{\elevenbb M}}
\def\R{\hbox{\elevenbb R}}
\def\S{\hbox{\elevenbb S}}
\def\T{\hbox{\elevenbb T}}
\def\Gr{Gr\"obner }
\def \bg #1 {\begin{tabular}{{#1}}}
\def \nd {\end{tabular}}
\newcommand \mwhile {{\bf while}\hspace{0.3cm}}
\newcommand \mfore {{\bf for\hspace{0.3cm}each}\hspace{0.3cm}}
\newcommand \mdo {{\bf do}\hspace{0.3cm}}

\newcommand \mif {{\bf if}\hspace{0.3cm}}
\newcommand \mthen {{\bf then}\hspace{0.3cm}}
\newcommand \melse {{\bf else}\hspace{0.3cm}}
\newcommand \mchoose {{\bf choose}\hspace{0.3cm}}

\newcommand \mbegin {{\bf begin}}
\newcommand \mend {{\bf end}}
\newcommand \bb {\hspace{0.3cm}}

\newcommand \hh {\hspace{1.0cm}}
\newcommand \hhh {\hspace{1.5cm}}
\newcommand \hhhh {\hspace{2.0cm}}
\newcommand \hhhhh {\hspace{2.5cm}}
\newcommand \hhhhhh {\hspace{3.0cm}}

\newtheorem{definition}{Definition}[section]
\newtheorem{corollary}[definition]{Corollary}
\newtheorem{proposition}[definition]{Proposition}
\newtheorem{example}[definition]{Example}
\newtheorem{theorem}[definition]{Theorem}

\begin{document}
\title{\bf Involutive Bases of Polynomial Ideals
}
\author{Vladimir P. Gerdt \\
       Laboratory of Computing Techniques and Automation\\
       Joint Institute for Nuclear Research\\
       141980 Dubna, Russia \\
       gerdt@jinr.dubna.su
\and
       Yuri A. Blinkov \\
       Department of Mathematics and Mechanics \\
       Saratov University \\
       410071 Saratov, Russia \\
       blinkov@scnit.saratov.su}
\date{}
\maketitle
\begin{abstract}In this paper we consider an algorithmic technique more
general than that proposed by Zharkov and Blinkov for the
involutive analysis of polynomial ideals. It is based on a new
concept of involutive monomial division which is defined for a
monomial set. Such a division provides for each monomial the
self-consistent separation of the whole set of variables into two
disjoint subsets. They are called multiplicative and
non-multiplicative. Given an admissible ordering, this separation
is applied to polynomials in terms of their leading monomials. As
special cases of the separation we consider those introduced by
Janet, Thomas and Pommaret for the purpose of algebraic analysis of
partial differential equations. Given involutive division, we
define an involutive reduction and an involutive normal form. Then
we introduce, in terms of the latter, the concept of involutivity
for polynomial systems. We prove that an involutive system is a
special, generally redundant, form of a Gr\"obner basis. An
algorithm for construction of involutive bases is proposed. It is
shown that involutive divisions satisfying certain conditions, for
example, those of Janet and Thomas, provide an algorithmic
construction of an involutive basis for any polynomial ideal. Some
optimization in computation of involutive bases is also analyzed.
In particular, we incorporate Buchberger's chain criterion to avoid
unnecessary reductions. The implementation for Pommaret division
has been done in Reduce.
\end{abstract}

\section{Introduction}
\noindent
In modern times the \Gr bases method invented by Buchberger~\cite{Buch65}
has become one of the most universal algorithmic
tools for analyzing and solving polynomial equations~\cite{Buch85,BWK93}.
Even in the general case, when the roots cannot be exactly computed, the
method is still able to obtain valuable information about the solutions. In
particular, it allows one to verify
compatibility of the initial equations and compute the dimension
of the solution space. For the last few
years notable progress has been achieved in extension
of the \Gr bases method to non-commutative~\cite{Mora,KRW90} and
differential algebras~\cite{Carra87,Ollivier}.

On the other hand, already by the early 20s the foundation of a constructive
approach to algebraic
analysis of partial differential equations was laid by
Riquier~\cite{Riquier} and Janet~\cite{Janet} giving, among other things,
answers to the
same
general questions of compatibility and dimension. Later on, this
approach, in the context of partial differential equations, was developed by
Thomas~\cite{Thomas} and
more recently by Pommaret~\cite{Pommaret78}. The main idea of the approach,
as with the computation of a \Gr basis, is rewriting the initial
differential
system into another, so-called, involutive form~\cite{Gerdt95}.

In the involutive approach, unlike the \Gr basis method,
independent variables for each equation are separated into two distinct
groups called multiplicative and non-multiplicative. Such a separation is
determined by the structure of the leading derivative terms. A differential
system is called involutive if its non-multiplicative derivatives are
algebraic consequences of multiplicative ones. In doing so,
Janet~\cite{Janet}, Thomas~\cite{Thomas} and Pommaret~\cite{Pommaret78}
used different separations of variables.

Zharkov and Blinkov~\cite{ZB93,ZB94} argued that the involutive
technique along with the \Gr bases one can be used in commutative
algebra. Based
on Pommaret definition of multiplicative and non-multiplicative
variables~\cite{Pommaret78}, they proved, among other things,
that an involutive basis is a \Gr one.
Moreover, their computational experience demonstrated a reasonably high
efficiency of the new algorithm when it terminates. The termination,
however, does not hold, generally, for positive dimensional ideals,
while for zero-dimensional ones it does for any degree-compatible monomial
orderings~\cite{ZB94}.
Apart from that, the Pommaret involutive form of Gr\"obner bases for
zero-dimensional polynomial ideals reveals a number of rather attractive
features~\cite{Zharkov94a}.

In the present paper we consider an algorithmic technique more
general than that proposed in~\cite{ZB93,ZB94} for the involutive
analysis of polynomial ideals. First of all, we introduce a new
concept of involutive monomial division (Sect.3) which leads to the
self-consistent separation of the whole set of variables into
multiplicative and non-multiplicative subsets. Given an admissible
ordering, the separation is applied to polynomials in terms of
their leading monomials. That concept generalizes
the particular
choice used by Janet~\cite{Janet}, Thomas~\cite{Thomas} and
Pommaret~\cite{Pommaret78} for analysis of partial differential
equations. We characterize also important properties of noetherity,
continuity and constructivity for involutive divisions (Sect.4).
Noetherity provides for the existence of a finite involutive basis
for any polynomial ideal. The other two properties allows one to
construct that basis algorithmically. It is shown that all the
above three divisions are continuous and constructive. Thomas and
Janet divisions are also noetherian whereas Pommaret division is
not.

Given an involutive division, we define an involutive
reduction and an involutive normal form (Sect.5). As this takes place, we
show
that much like the Pommaret normal form, investigated in~\cite{ZB93},
the general involutive normal form is also unique and linear.
Then we define involutive systems by analogy with
differential equations (Sect.6). To be involutive, systems are required to
satisfy the involutivity conditions, which form the basis for
their algorithmic construction.

We prove (Sect.7) that any involutive basis, if it exists, is a
special, generally extended, form of the reduced Gr\"obner basis.
Though it is unique for Pommaret division~\cite{ZB94}, generally,
it may not be the case, as it is shown by an explicit example. We
propose an algorithm for construction of involutive polynomial
bases (Sect.8). Its correctness is proved for any continuous
involutive division and for arbitrary admissible monomial ordering,
while its termination holds, generally, for noetherian divisions.
The algorithm is an improved and generalized version of one
proposed in~\cite{ZB94,Zharkov94a}, and has been implemented in
Reduce for Pommaret division. The main improvement is the
incorporation of Buchberger's chain criterion~\cite{Buch79}.

\section{Preliminaries}
\noindent
Let $\R=K[x_1,\ldots,x_n]$ be a polynomial ring over the field $K$ of
characteristic zero. In this paper we use the notations:
\vskip 0.2cm
\noindent
\hspace*{0.5cm}$f,g,h,p,q$\hspace*{0.25cm} are polynomials in $\R$. \\
\hspace*{0.5cm}$a,b,c$\hspace*{0.25cm} are elements in $K$.\\
\hspace*{0.5cm}$F,G,H$\hspace*{0.25cm} are finite subsets of $\R$.\\
\hspace*{0.5cm}$\N$\hspace*{0.25cm} is the set of non-negative integers.\\
\hspace*{0.5cm}$\M=\{\ x_1^{d_1}\cdots x_n^{d_n}\ |\
d_i\in \N,\,i=1,\ldots,n\ \}$\hspace*{0.25cm}
is the set of monomials in $\R$.\\
\hspace*{0.5cm}$\T=\{\ a\,u\ |\   u\in \M,\,a\in K\ \}$\hspace*{0.25cm} is the
set of terms in $\R$.\\
\hspace*{0.5cm}$u,v,w,s,t$ are monomials or terms with nonzero coefficients.\\
\hspace*{0.5cm}$U,V,W$\hspace*{0.25cm} are finite subsets of $\M$.\\
\hspace*{0.5cm}$deg_i(u)$\hspace*{0.25cm} is the degree of $x_i$ in $u$.\\
\hspace*{0.5cm}$deg(u)$\hspace*{0.25cm} is the total degree of $u$.\\
\hspace*{0.5cm}$cf(f,u)\in K$\hspace*{0.25cm} is the coefficient of the term $u$
 of the polynomial $f$.\\
\hspace*{0.5cm}$Id(F)$\hspace*{0.25cm} is the ideal in $\R$ generated by the
 polynomial set $F$.\\
\hspace*{0.5cm}$\succ$\hspace*{0.25cm} is an admissible monomial ordering with
$x_1\succ x_2\succ\cdots\succ x_n$. \\
\hspace*{0.5cm}$lt(f)$\hspace*{0.25cm} is the leading term of $f$ w.r.t. the
ordering $\succ$.\\
\hspace*{0.5cm}$lc(f)=cf(f,lt(f))$\hspace*{0.25cm} is the leading
coefficient of $f$.\\
\hspace*{0.5cm}$lm(f)=lt(f)/lc(f)$\hspace*{0.25cm} is the leading monomial of $f$.\\
\hspace*{0.5cm}$lm(F)=\{\ lm(f)\ |\   f\in F\ \}$\hspace*{0.25cm} is the set of
the leading monomials of $F$.\\
\hspace*{0.5cm}$lcm(F)$\hspace*{0.25cm} is the least common multiple of the
set $\{\ lm(f)\ |\   f\in F\ \}$.
\vskip 0.2cm
\noindent
If the monomial $u$ divides the monomial $v$ we shall write $u|v$.

\section{Involutive Monomial Division}
\noindent
\begin{definition}
{\em We shall say that an {\em involutive  division} $L$
or
 $L-${\em division}
 is given on $\M$ if for any finite set $U\subset \M$ a
 relation $|_L$ is defined on $U\times \M$
 such that for any $u,u_1\in U$ and $v,w\in \M$ the following holds:
\begin{enumerate}
\renewcommand{\theenumi}{(\roman{enumi})}
\item $u|_L  w$ implies $u|  w$.
\item $u|_L  u$ for any $u\in U$.
\item $u|_L (uv)$ and $u|_L (uw)$ if and only if $u|_L (uvw)$.
\item If $u|_L w$ and $u_1|_L w$, then $u|_L u_1$ or
  $u_1|_L u$.
\item If $u|_L u_1$ and $u_1|_L  w$, then $u|_L
 w$.
\item If $V\subseteq U$ and $u\in V$, then $u|_L w$ w.r.t. $U$ implies
 $u|_L w$ w.r.t. $V$.
\end{enumerate}
If $u|_L (w=uv)$, we say $u$ is an {\em involutive divisor} of $w$,
$w$ is an {\em involutive multiple} of $u$, and $v$ is {\em  multiplicative}
for $u$. In such an event we shall write $w=u\times v$. If $u$ is a
conventional divisor of $w$ but not an involutive one we shall write,
as usual, $w=u\cdot v$. Then $v$ is said to be {\em non-multiplicative}
for $u$.
} \label{inv_div}
\end{definition}
The conventional monomial division, obviously, satisfies condition (iv) only
in the univariate case. The simplest bivariate example: $x|(xy)$
and $y|(xy)$ but $\neg x|y$ and $\neg y|x$.

Definition~\ref{inv_div} for each $u\in U$ provides separation of the set
of variables $$\{x_1,\ldots,x_n\}=M_L(u,U)\cup NM_L(u,U)$$
into two disjoined subsets
$(M_L(u,U)\cap NM_L(u,U)=\emptyset)$
of {\em multiplicative} $M_L(u,U)$ and
{\em non-multiplica\-ti\-ve} $NM_L(u,U)$ variables. It is convenient to define
an involutive division for a monomial set just by specifying the
subsets of multiplicative and non-multiplicative variables to satisfy the
conditions (iv)-(vi). The other
conditions will be fulfilled by the construction.

Given an involutive division $L$ and a finite set $U$, for each $u\in U$
let $L(u,U)\subseteq \M$ be the set of multiplicative monomials for $u$,
that is,
\begin{equation}
u|_L v\ \Longleftrightarrow \ v\in uL(u,U). \label{mult_set}
\end{equation}
Then it is easy to see that Definition~\ref{inv_div} admits another
form:
\begin{definition}
{\em An {\em involutive division} $L$ on $\M$ is given, if for any finite
$U\subset \M$ and for any $u\in U$ there is given a
submonoid $L(u,U)$ of
$\M$ satisfying the conditions:
\begin{enumerate}
\renewcommand{\theenumi}{(\alph{enumi})}
\item If $w\in L(u,U)$ and $v|w$, then $v\in L(u,U)$.
\item If $u,v\in U$ and $uL(u,U)\cap vL(v,U)\not=\emptyset$, then
 $u\in vL(v,U)$ or $v\in uL(u,U)$.
\item If $v\in U$ and $v\in uL(u,U)$, then $L(v,U)\subseteq L(u,U)$.
\item If $V\subseteq U$, then $L(u,U)\subseteq L(u,V)$ for all $u\in V$.
\end{enumerate}
} \label{3.1a}
\end{definition}

We consider three different examples of involutive division
introduced by Janet~\cite{Janet}, Thomas~\cite{Thomas} and
Pommaret~\cite{Pommaret78} for
analysis of algebraic differential equations. In doing so, we give,
firstly, the definition of multiplicative and non-multiplicative
variables for each of the divisions, and, secondly, prove the
fulfillment of the three extra conditions (iv)-(vi) in
Definition~\ref{inv_div} equivalent to (b)-(d) in Definition~\ref{3.1a}.

\begin{definition}{\em{\em Thomas division}~\cite{Thomas}.
Given a finite set $U$, let
$$h_i(U)=max\{\ deg_i(u)\ |\   u\in U\ \}\,.$$
A variable $x_i$ is considered as multiplicative for $u\in U$ if
$deg_i(u)=h_i(U)$ and non-multiplicative, otherwise. }
\label{Thomas}
\end{definition}

\begin{definition}{\em {\em Janet division}~\cite{Janet}. Let $U$
be a finite set. For each $1\leq i\leq n$ divide $U$ into groups
labeled by non-negative integers $d_1,\ldots,d_i$:
$$[d_1,\ldots,d_i]=\{\ u\ \in U\ |\ deg_j(u)=d_j,\ 1\leq j\leq i\ \}.$$
A variable $x_i$ is multiplicative for $u\in U$ if
$i=1$ and $deg_1(u)=max\{\ deg_1(v)\ |\ v\in U\ \}$,
or if $i>1$, $u\in [d_1,\ldots,d_{i-1}]$ and
$$ deg_i(u)=max\{\ deg_i(v)\ |\ v\in [d_1,\ldots,d_{i-1}]\ \}.$$
} \label{Janet}
\end{definition}

\begin{definition}{\em
{\em Pommaret division}~\cite{Pommaret78}. For a monomial
$x_1^{d_1}\cdots x_k^{d_k}$ with $d_k>0$ the variables $x_j$ with $j\geq k$
are considered as multiplicative and $x_j$ with $j<k$ as non-multiplicative.
For $u=1$ all the variables are multiplicative.}
\label{Pommaret}
\end{definition}
We note that
\begin{itemize}
\item Thomas division does not depend on the ordering on the variables $x_i$.
Janet and Pommaret divisions, as defined, are based on the ordering
of the variables assumed in Sect.2.
\vskip 0.2cm
\item The separation of variables into multiplicative and non-multiplicative
ones for Thomas and Janet
divisions are defined in terms of the whole set $U$.
Contrastingly, Pommaret division is determined in terms of the monomial
itself, regardless of the others, and, by this reason, admits extension
to infinite monomial sets, unlike Thomas and Janet divisions.
\end{itemize}
To distinguish the above divisions the related subscripts $T,J,P$ will be
used.

\begin{proposition}
Thomas, Janet and Pommaret monomial divisions are involutive.
\end{proposition}

\vskip 0.2cm
\noindent
{\bf Proof}\ \
According to the above remark we must prove that the conditions
(iv)-(vi) in Definition~\ref{inv_div} are satisfied.

Let $u\in U$ be a Thomas divisor of $w\in \M$, that is, $w=u\times
v$. Then $deg_i(v)=deg_i(w)-h_i(U)$ if $deg_i(w)\geq h_i(U)$ and
$deg_i(v)=0$ if $deg_i(w)< h_i(U)$. Thus, if $w$ has an involutive
divisor $u$, then $w/u$ is uniquely defined, and, hence, $u$ is
unique in $U$. It implies also the property (v) for Thomas
division, since $u|_{T}v$ for $u,v\in U$ if and only if $u=v$. The
property (vi) also holds since any $h_i$ for $V$ is less than or
equal to the corresponding $h_i$ for $U$.

Let now $u,v\in U$ be two different Janet divisors of $w$, such that
$deg_i(u)=deg_i(v)=d_i$ for $1\leq i< k\leq n$ and assume, for definiteness,
that $deg_k(u)>deg_k(v)$. Then, since both $u,v$ are members of the same
group
$[d_1,\ldots,d_{k-1}]$, the variable $x_k$ is non-multiplicative for $v$.
Hence, if $u$ is a Janet divisor of $w$ such that
$deg_k(w)\geq deg_k(u)>deg_k(v)$, then
$v$ is not Janet divisor of $w$. In other words, similar to Thomas division,
any monomial $w\in \M$ cannot have different Janet divisors in any set $U$.
A monomial group may only be decreased by diminishing the set $U$, which
implies the relation (vi).

Lastly, consider a Pommaret divisor $u$ of the monomial
$w=x_1^{d_1}\cdots x_m^{d_m}$ with $m\leq n$. By definition, $u$
constitutes a left subset of the string representation for $w$ as it is shown.
\begin{equation}
w=\overbrace{\underbrace{x_1\cdots x_1}_{d_1}\cdots}^u\cdots
\underbrace{x_m\cdots x_m}_{d_m}\,. \label{string}
\end{equation}
It makes evident the fulfillment of the conditions (iv) and (v) for Pommaret
division while the condition (vi) trivially holds since the division
does not depend on the set $U$ at all.
\hfill{\Box}

\begin{proposition} For any finite set $U$ and for any
$u\in U$, the inclusion $M_T(u,U)\subseteq M_J(u,U)$ and, respectively,
$NM_J(u,U)\subseteq NM_T(u,U)$ holds.
\label{id_T_J}
\end{proposition}

\vskip 0.2cm
\noindent
{\bf Proof}\ \  If $x_i\in M_J(u,U)$, $u\in [d_1,\ldots,d_{i-1}]$, then,
by definition,
$$deg_i(u)=max\{\ deg_i(v)\ |\   v\in [d_i,\ldots,d_{i-1}]\ \}\leq
max\{\ deg_i(v)\ |\   v\in U\ \}\,.$$
Hence, $x_i\in M_T(u,U)$ implies $x_i\in M_J(u,U)$.
\hfill{\Box}

\begin{definition}{\em A set $U$ is called {\em involutively
autoreduced} with respect to division $L$ or {\em $L-$autoreduced}
if it does not contain elements $L-$divisible by other
elements in $U$.
}\end{definition}

\begin{proposition}If $U$ is $L-$autoreduced, then any
 monomial $w\in \M$ has at most one $L-$involutive divisor in $U$. \label{unique_id}
\end{proposition}

\vskip 0.2cm
\noindent
{\bf Proof}\ \  This follows immediately from the property (iv) of involutive
division. In terms of Definition~\ref{3.1a} it means that
$uL(u,U)\cap vL(v,U)=\emptyset$ for all distinct $u,v\in U$, if
$U$ is involutively autoreduced.
\hfill{\Box}

\begin{proposition}~\cite{Zharkov94b}. If a set $U$ is autoreduced
with respect to Pommaret division, then for any
$u\in U$ $M_P(u,U)\subseteq M_J(u,U)$ and
$NM_J(u,U)\subseteq NM_P(u,U)$, respectively.  \label{id_P_J}
\end{proposition}

\vskip 0.2cm
\noindent
{\bf Proof}\ \ Let $u=x_1^{d_1}\cdots x_k^{d_k}\in \M$ be a monomial with
$d_k>0$
and $v\in U$ be its Pommaret divisor. Then, as follows from the
representation (\ref{string}),
$v=x_1^{d_1}\cdots x_{m-1}^{d_{m-1}}x_m^r$ with $1\leq m\leq k$ and
$1\leq r\leq d_m$. It means that
$v\in [d_1,\ldots,d_{m-1}]$. Since $U$ is autoreduced by Pommaret division,
there are no other members of the same group with degree in
$x_m$ higher than $r$. Therefore, $v$ is also a Janet divisor of $u$,
and $u/v$, being Pommaret multiplicative for $v$, is also Janet
multiplicative.
\hfill{\Box}

\begin{example} {\em $U=\{xy,y^2,z\}$ ($x\succ y\succ z$).
\vskip 0.3cm
\begin{center}
\bg {|c|c|c|c|c|c|c|} \hline\hline
 & \multicolumn{2}{c|}{Thomas} &\multicolumn{2}{c|}{Janet}
& \multicolumn{2}{c|}{Pommaret} \\ \cline{2-7}
monomial & $M_T$ & $NM_T$  &   $M_J$   & $NM_J$   & $M_P$   & $NM_P$ \\ \hline
 $xy$  & $x$ & $y,z$ & $x,y,z$ & $-$   & $y,z$  &  $x$    \\
 $y^2$ & $y$ & $x,z$ & $y,z$   & $x$    & $y,z$ &  $x$   \\
 $z$   & $z$ & $x,y$ &  $z$    & $x,y$  & $z$   & $x,y$  \\ \hline \hline
\nd
\end{center} \label{exm_1}
}
\end{example}

\section{Involutive Monomial Sets}

\begin{definition}{\em Given an involutive division $L$, a set $U$ is
 called {\em involutive} with respect to $L$ or {\em $L-$involutive}, if
 any multiple of some element $u\in U$, is also ($L-$)involutively
 multiple of an element $v\in U$, generally, different from $u$. It means
 that
\begin{equation}
 (\forall u\in U)\ (\forall w\in \M)\ (\exists v\in U)\
\ [\ v|_L(uw)\ ]\, \label{imset}
\end{equation}
or, in accordance with~(\ref{mult_set}) and Definition~\ref{3.1a},
$$\cup_{u\in U}\,u\,\M = \cup_{u\in U}\,u\,L(u,U).$$
} \label{inv_mset}
\end{definition}

\begin{definition}{\em We shall call the set $\cup_{u\in U}\,u\,\M$ the {\em
cone} generated by $U$ and denote it by $C(U)$. The set
$\cup_{u\in U}\,u\,L(u,U)$ will be called the {\em involutive cone} of $U$
with respect to $L$ and denoted by $C_L(U)$.
} \label{cone}
\end{definition}
Thus, the set $U$ is $L-$involutive if and only if its cone $C(U)$
coincides with its involutive cone $C_L(U)$.

\begin{definition}{\em A finite $L-$involutive set $\tilde{U}\subset \M$
 will be called $L-${\em completion} of a set $U\subseteq \tilde{U}$ if
 $C(\tilde{U})=C(U)$. If there exists an $L-$completion $\tilde{U}$ of
 the set $U$,
 then the latter is said to be {\em finitely generated} with respect
 to $L$. An involutive division $L$ is called {\em noetherian} if
 every finite set $U$ is finitely generated.
} \label{id_noetherian}
\end{definition}

\begin{proposition} Given a noetherian involutive division $L$, every
 monomial ideal $U$ has a finite involutive basis.
\label{cr_finite_imb}
\end{proposition}

\vskip 0.2cm
\noindent
{\bf Proof}\ \
 This is an immediately consequence of Definition~\ref{id_noetherian}
 and Dickson's lemma~\cite{BWK93}.
\hfill{\Box}

\begin{proposition} Thomas and Janet divisions are noetherian.
\label{pr_finite}
\end{proposition}

\vskip 0.2cm
\noindent
{\bf Proof}\ \  Given a finite set $U$, consider the monomial
$h=x_1^{h_1}\cdots x_n^{h_n}$ where, as given in the definition of Thomas
division, $h_i=max\{\ deg_i(u)\ |\   u\in U\ \}$,
and form the finite set $V\subset \M$ of all the different monomials
$v$ such that $v|h$ and $u|v$ for some $u\in U$. The set $V$, which
contains, in particular, the monomial $h$ and the initial set $U$, is
involutive for Thomas division. Indeed, let $w=x_1^{d_1}\cdots x_n^{d_n}$
be a multiple of some $u\in V$. If $w\in V$, then, obviously, $w\in C_T(V)$.
Otherwise, let $\{d_{i_1},\ldots,d_{i_k}\}$ ($k\leq n$) be the nonempty
set which contains all the exponents $d_i$ $(1\leq i\leq n)$ in $w$ such
that $d_{i_1}>h_{i_1},\ldots,d_{i_k}>h_{i_k}$. Then there exists $v\in V$
satisfying
$$w=v\, x_{i_1}^{d_{i_1}-h_{i_1}}\cdots x_{i_k}^{d_{i_k}-h_{i_k}}\,.$$
Since $deg_{i_1}(v)=h_{i_1},\ldots,deg_{i_k}(v)=h_{i_k}$, $v$ is a
Thomas involutive divisor of $w$, and, hence, $w\in C_T(V)$.

Furthermore, from Proposition~\ref{id_T_J} it follows that there is
a set of $V_1\subseteq V$ which is a Janet completion of $U$.
\hfill{\Box}

\begin{definition}{\em  Multiplication of a monomial $u\in U$ by a
 variable $x$ is called a {\em prolongation} of $u$. Given an involutive
 division specified by the set $U$, the prolongation
 is called {\em multiplicative} if $x$ is multiplicative for $u$ and
{\em non-multiplicative}, otherwise.}
\label{m_prolongation}
\end{definition}
In the construction of involutive sets the following concept of local
involutivity plays the crucial role and admits the direct extension to
polynomial sets (see Sect.6).

\begin{definition}
 {\em A set $U$ is called {\em locally involutive} with respect
 to the involutive division $L$ if any non-multiplicative prolongation
 of any element in $U$ has an involutive divisor in $U$, that is,
\begin{equation}
(\forall u\in U)\ (\forall x_i\in NM_L(u,U))\ (\exists v\in U)\
\ [\ v|_L(u\cdot x_i)\ ]\, \label{mon_inv_cond}
\end{equation}
}
\label{mon_loc_inv}
\end{definition}
In accordance with Definition~\ref{inv_mset},
the conditions~(\ref{mon_inv_cond}), apparently, are necessary for
involutivity of $U$. Generally, however, they are not sufficient, as the
next simple example shows.

\begin{example}{\em Let $L$ be an
involutive division on $\M\subset K[x,y,z]$ defined by the table
\vskip 0.3cm
\begin{center}
\bg {|c|c|c|} \hline\hline
monomial         & $M$   & $NM$  \\ \hline
 $1$             & $x,y,z$ & $-$   \\
 $x$             & $x,z$ & $y$   \\
 $y$             & $x,y$ & $z$   \\
 $z$             & $y,z$ & $x$   \\
 $u\in \M\ |\ deg(u)\geq 2$ & $-$ & $x,y,z$ \\ \hline \hline
\nd
\end{center}
\vskip 0.3cm
It is easy to see that all properties listed in
Definition~\ref{inv_div}~(\ref{3.1a}) are satisfied, and the set
$U=\{x,y,z\}$ is locally involutive. For instance, $x\cdot y=y\times x$.
However, $U$ is not involutive since none $u\in \M$ with
$deg_x(u)>0,deg_y(u)>0,deg_z{u}>0$, e.g. $xyz$, has involutive divisors
in $U$.
} \label{c_example}
\end{example}
The following definition and theorem enable one to reveal
involutive divisions providing involutivity of every locally
involutive set.

\begin{definition}
{\em An involutive division $L$ will be called {\em continuous} if
for any finite set $U$ and for any
finite sequence $\{u_i\}_{(1\leq i\leq k)}$ of elements
in $U$ such that
\begin{equation}
(\forall \,i< k)\ (\exists x_j\in NM_L(u_i,U))\ \ [\ u_{i+1}|_L u_i\cdot
x_j\ ] \label{cont_cond}
\end{equation}
the inequality $u_i\neq u_j$ for $i\neq j$ holds.
} \label{def_cont}
\end{definition}

\begin{theorem}
If an involutive division $L$ is continuous then
local involutivity of any set $U$ implies
its involutivity.
\label{th_cont}
\end{theorem}

\vskip 0.2cm
\noindent
{\bf Proof}\ \
Let set $U$ be locally involutive, and such that any
sequence in
$U$ satisfying~(\ref{cont_cond}) has no coinciding elements.
We must prove that $U$ satisfies~(\ref{imset}).
Take any $u\in U$ and any $w\in \M$ and show that there is $v\in U$
such that $v|_L (uw)$.
If $u|_L(uw)$ we are done. Otherwise, there is $x_{k_1}\in NM_L(u,U)$
such that $w$ contains $x_{k_1}$. Then $u\cdot x_{k_1}$ has an involutive
divisor $v_1\in U$. If $v_1|_L(uw)$ we are
done. Otherwise, there are $x_{k_2}\in NM_L(v_1,U)$ and $v_2\in U$
such that $uw/v_1$ contains $x_{k_2}$ and $v_2|_L (v_1\cdot x_{k_2})$.
Going on, we obtain the sequence $u,v_1,v_2,\ldots $ of elements in $U$
satisfying~(\ref{cont_cond}). By construction, each element of
the sequence divides $uw$. Since all the elements are distinct and
$uw$ has a finite number of distinct divisors, it follows that the above
sequence in $U$ is finite, and, hence, it ends up with an involutive divisor
of $uw$.
\hfill{\Box}

\begin{corollary} Thomas, Janet and Pommaret divisions are
continuous.
\label{th_cont_TJP}
\end{corollary}

\vskip 0.2cm
\noindent
{\bf Proof}\ \
Let $U$ be a finite set, and $\{u_i\}_{(1\leq i\leq k)}$ be a sequence
of elements in $U$ satisfying the conditions~(\ref{cont_cond}).
We shall show that there cannot be coinciding elements in the sequence for
three divisions.

It is ease to see that $u_{i+1}|_T (u_i\cdot x_{k_i})$ implies
$u_{i+1}=u_i\cdot x_{k_i}$. Indeed,
suppose that $u_i\cdot x_{k_i}=u_{i+1}\times v_{i+1}$ what means
$\neg x_{k_i}|v_{i+1}$.
If $v_{i+1}$ would contain any other variable $x_{j_i}$,
then it would mean that $deg_{x_{j_i}}(u_i)>deg_{x_{j_i}}(u_{i+1})$,
and, hence, $x_{j_i}$ could not be multiplicative for $u_{i+1}$.
Therefore, any Thomas sequence satisfying~(\ref{cont_cond}) consists
of distinct elements.

If $u_{i+1}|_J (u_i\cdot x_{k_i})$, then from definition of Janet division
it follows that $u_{i+1}\succ_{Lex} u_i$, where $\succ_{Lex}$
is the lexicographical ordering corresponding to the choice of
variable order $x_1\succ x_2\succ  \cdots \succ x_n$ as assumed in
Sect.2. It is now
obvious that $u_i\neq u_j$ for $i\neq j$ for Janet division.

Let now $u_{i+1}|_P(u_i\cdot x_{k_i})$. Then
the representation~(\ref{string}) shows clearly that
$u_{i+1}\succ_{RevLex}u_i$ where $\succ_{RevLex}$ is the reverse
lexicographical ordering on $\M$ induced by the assumed variable order.
\hfill{\Box}

\vskip 0.2cm
\noindent
With an eye to the below described algorithms
based on examination of non-multiplica\-tive prolongations only, we
impose, in addition to continuity, one more requirement on an
involutive division.

\begin{definition}
{\em  We shall say that a continuous involutive division $L$
 is {\em constructive} if for any $U\subset \M$, $u\in U$,
 $x_i\in NM_L(u,U)$ such that $u\cdot x_i\not \in C_L(U)$ and
\begin{equation}
 (\forall v\in U)\ (\forall x_j\in NM_L(v,U))\ (v\cdot x_j | u\cdot x_i,\
 v\cdot x_j\neq u\cdot x_i)\ \ [\ v\cdot x_j\in C_L(U))\ ]
  \label{constr1}
\end{equation}
the following condition holds:
\begin{equation}
 (\forall w\in C_L(U))
 \ \ [\ u\cdot x_i\not \in wL(w,U\cup \{w\})\ ].
\label{constr}
\end{equation}
} \label{def_constr}
\end{definition}

\begin{proposition} Thomas, Janet and Pommaret divisions are
constructive.
\label{pr_constr_TJP}
\end{proposition}

\vskip 0.2cm
\noindent
{\bf Proof}\ \  Let $T$ be Thomas division. Suppose there is $u_1\in U$,
 and $v\in T(u_1,U)$ such that
 $u\cdot x_i = u_1v\times w$, $w\in T(u_1v,U\cup \{u_1v\})$.
 From Definition~\ref{Thomas} it
 follows that if there exists $x_j|w$ and $\neg x_j|v$ for some
 $1\leq j\leq n$, then $x_j\in M_T(u_1,U)$. This implies
 $w\in T(u_1,U)$ and $u\cdot x_i\in u_1T(u_1,U)$.

 Consider now Janet division $J$, and let $u\cdot x_i$ be a
 non-multiplicative
 prolongation which has no Janet divisors in $U$, and for which the
 condition~(\ref{constr1}) holds. Assume for a contradiction
 that there is $u_1\in U$ and $v\in J(u_1,U)$
 satisfying
  $$u\cdot x_i = u_1v\times w_1,\quad w_1\in J(u_1v,U\cup \{u_1v\}).$$
 Because $v\neq 1$ and $w_1\neq 1$, select
 minimal $j,m$ such that $x_j|v$ and $x_m|w_1$. It is easy to see that
 $i < min\{j,m\}$. Otherwise, by Definition~\ref{Janet}, we would have
 either $x_j\not \in J(u_1,U)$ if $i\geq j$ or $x_m\not \in
 J(u_1v,U\cup \{u_1v\})$ if $i>m$. Note that the equality $i=m$
 impossible since $U\cup \{u_1v\}$ as well as any other monomial set is
 Janet autoreduced. Thus, $u_1\succ_{Lex} u$ where $\succ_{Lex}$ is
 the lexicographical ordering
 induced by the variable order $x_1\succ \cdots \succ x_n$. If monomial
 $ux_i$ is obtained by non-multiplicative prolongations of several elements
 in $U$, then we suppose that $u$ is lexicographically maximal from
 all of them. Since $w_1$ is non-multiplicative for $u_1$, there is
 $x_{k_1}|w_1$ such that $x_{k_1}\in NM_J(u_1,U\cup \{u_1v\})$. Then,
 by condition (\ref{constr1}), we can rewrite
 $$ u\cdot x_i=(u_1\cdot x_{k_1})\frac{vw_1}{x_{k_1}}=
 (u_2\times w_2)\frac{vw_1}{x_{k_1}}=(u_2\cdot x_{k_2})
 \frac{vw_1w_2}{x_{k_1}x_{k_2}}=\cdots\,, $$
 where $u\prec_{Lex}u_1\prec_{Lex}u_2\prec_{Lex}\cdots$. Continuity of
 Janet division implies termination of this chain with some $u_l\in U$ such that
 $u\cdot x_i\in u_lJ(u_l,U)$ what contradicts our assumption
 $u\cdot x_i\not \in C_J(U)$.

 For Pommaret division condition~(\ref{constr}) follows directly
 from the property (v) in Definition~\ref{inv_div}.

\hfill{\Box}

\begin{theorem} Let $U$ be a non-involutive finitely generated set with
 respect to a constructive division $L$. Then there is a
 procedure of completing $U$ to an $L-$involutive set $\tilde{U}\supset U$
 based on enlargement of $U$ by non-multiplicative prolongations of its
 elements.
\label{th_inv_completion}
\end{theorem}

\vskip 0.2cm
\noindent
{\bf Proof}\ \
 Given $U$, by Definition~\ref{id_noetherian}, there exists a finite
 $L-$completion $\tilde{U}$ of $U$.
 We claim that $\tilde{U}$ contains some non-multiplicative prolongations
 of elements in $U$. Assume for
 a contradiction that there are no such elements in $\tilde{U}$.
 Since set $U$ is not involutive, there exist non-multiplicative
 prolongations of elements in $U$ which have no $L-$divisors in $U$.

 Take any admissible ordering $\prec$ and select
 $u\in U$ with a non-multiplicative prolongation
 $u\cdot x_i$ which is not $L-$multiple of any element
 in $U$, and which is the lowest with respect to
 $\prec$. Because $\tilde{U}$ is involutive, and, by the above
 assumption, $u\cdot x_i\not \in \tilde{U}$, there is
 $v\in \tilde{U}\setminus U$ and $1\prec w\in \M$ such that
 $u\cdot x_i=v\times w$, $w\in L(v,\tilde{U})$. From the condition
 $C(U)=C_L(\tilde{U})$ it follows that $v$ is multiple of some
 $u_1\in U$ with $deg(u_1)<deg(v)$.

 Show that $v\in C_L(U)$. If $u_1$ $L-$divides $v$, then we are done.
 Otherwise, there exists $x_{k_1}|(v/u_1)$, $x_{k_1}\in NM_L(u_1,U)$,
 and we can rewrite
 $$v=u_1\cdot \frac{v}{u_1}=(u_1\cdot x_{k_1})\frac{v}{u_1x_{k_1}}=
 (u_2\times w_2)\frac{v}{u_1x_{k_1}}=
 (u_2\cdot x_{k_2})\frac{vw_2}{u_1x_{k_1}x_{k_2}}=\cdots $$
 until, by continuity of $L$, we come to an involutive divisor $u_m\in U$ of
 $v$ at some step of this rewriting procedure.  This contradicts
 the constructivity condition~(\ref{constr}), and, hence
 $u\cdot x_i\in \tilde{U}$.

 Now instead of $U$ take $U_1=U\cup\{u\cdot x_i\}$ where
 $u_1\in U$ and $u_1\cdot x_{i_1}\in \tilde{U}\setminus U$.
 If set $U_1$ is not involutive, then it can be further completed by
 the corresponding lowest non-multiplicative prolongation in $U_1$.
 Since the set $\tilde{U}$ is finite, by repeating this completion
 procedure, in a finite number of steps we construct the set
 $\bar{U}\subseteq \tilde{U}$ which is an $L-$completion of $U$.
\hfill{\Box}

\vskip 0.2cm
\noindent
 As an immediate consequence of the above described procedure
 of completing a set $U$ by non-multiplicative prolongations of its
 elements we have the following corollary.

\begin{corollary}
 If $U$ is a finitely generated set with respect to a constructive
 involutive division, then there is the unique minimal involutive
 completion $\bar{U}$ of $U$ such that for any other completion
 $\tilde{U}$ the inclusion $\bar{U}\subseteq \tilde{U}$ holds.
\end{corollary}
The following algorithm, given a constructive division $L$, computes
the minimal involutive completion $\tilde{U}$ for any finitely
generated set $U$ and any fixed admissible ordering $\prec$. Its
{\em correctness} and {\em termination} are provided
by Theorem~\ref{th_inv_completion}.

\vskip 0.3cm
\noindent
\hh Algorithm {\bf InvolutiveCompletion:}
\vskip 0.2cm
\noindent
\hh {\bf Input:}  $U$, a finite monomial set
\vskip 0.0cm \noindent
\hh {\bf Output:} $\tilde{U}$, an involutive completion of $U$
\vskip 0.0cm \noindent
\hh \mbegin
\vskip 0.0cm \noindent
\hhh $\tilde{U}:=U$
\vskip 0.0cm \noindent
\hhh \mwhile exist $u\in \tilde{U}$ and $x\in NM_L(u,\tilde{U})$ such that
\vskip 0.0cm \noindent
\hhhh $u\cdot x$ has no involutive divisors in $\tilde{U}$\bb \mdo
\vskip 0.0cm \noindent
\hhhh \mchoose such $u$ and $x$ with the lowest $u\cdot x$ w.r.t. $\prec$
\vskip 0.0cm \noindent
\hhhh $\tilde{U}:=\tilde{U}\cup \{u\cdot x\}$
\vskip 0.0cm \noindent
\hhh \mend
\vskip 0.0cm \noindent
\hh \mend
\vskip 0.3cm
\noindent

\begin{example} {\em (Continuation of Example~\ref{exm_1}). The
 minimal involutive bases of the set $U=(xy,y^2,z)$ ($x\succ y\succ z$)
 for Thomas, Janet and Pommaret divisions are
\begin{eqnarray*}
&\bar{U}_T&=\{xy,y^2,z,xz,yz,xy^2,xyz,y^2z,xy^2z\}\,,\\
&\bar{U}_J&=\{xy,y^2,z,xz,yz\}\,,\\
&\bar{U}_P&=\{xy,y^2,z,xz,yz,x^2y,x^2z,\ldots,x^ky,\ldots,x^mz,\ldots\}\,,
\end{eqnarray*}
\vskip 0.2cm
\noindent
where $k,m\in \N$. These bases can be easily derived from $U$ using
algorithm {\bf InvolutiveCompletion}. Note that $\bar{U}_J\subset
\bar{U}_T$ and $\bar{U}_J\subset \bar{U}_P$ in agreement with
Propositions~\ref{id_T_J} and~\ref{id_P_J}. This example explicitly
shows that Pommaret division is not noetherian. However, for
another ordering $z\succ y\succ x$ the set $U$ is finitely
generated, and then $\bar{U}_P=U$.} \label{exm_2}
\end{example}

\section{Polynomial Reduction}
\noindent
In this section we generalize the results obtained in~\cite{ZB93,ZB94}
for Pommaret division to arbitrary involutive division.

\begin{definition}{\em Given a finite polynomial set $F\subset\R$
and an admissible ordering $\succ$, the concept of multiplicative and
non-multiplicative variables for $f\in F$ is to be defined in terms of
$lm(f)$ and the leading monomial set $lm(F)$.}
\end{definition}
Therefore, as soon as we have polynomials rather than monomials, any involutive
division is to be determined on the basis of some admissible ordering, even
when it does not depend on the latter for the pure monomial case, as with
Thomas division.

The concepts of involutive polynomial reduction and involutive normal form
are introduced similar to their conventional analogues (Buchberger, 1985) with
the use of involutive division instead of the conventional one.

\begin{definition} {\em Let $L$ be an involutive division $L$ on $\M$, and
 let $F$ be a finite set of polynomials. Then we shall say:
 \begin{enumerate}
 \renewcommand{\theenumi}{(\roman{enumi})}
 \item $p$ is {\em $L-$reducible} {\em modulo} $f\in F$ if
  $p$ has a term $t=a\,u\in \T$ ($a\neq 0$) such that $u=lm(f)\times v$,
  $v\in L(lm(f),lm(F))$.
  It yields the {\em $L-$reduction} $p\rightarrow g=p-(a/lc(f))\,f\times v$.
 \item $p$ is {\em $L-$reducible modulo} $F$ if there exists $f\in F$ such
  that $p$ is $L-$reducible modulo $f$.
 \item $p$ is {\em in $L-$normal form modulo $F$} if
  $p$ is not $L-$reducible modulo $F$.
\end{enumerate}
} \label{inv_red}
\end{definition}
We denote an $L-$ normal form of $p$ modulo $F$ by $NF_L(p,F)$.
In contrast, a conventional normal form will be denoted by $NF(p,F)$.
As an involutive normal form algorithm one can use, for example, the
following:

\vskip 0.3cm
\noindent
\hh Algorithm {\bf InvolutiveNormalForm:}
\vskip 0.2cm
\noindent
\hh {\bf Input:}  $p,\,F$
\vskip 0.0cm \noindent
\hh {\bf Output:} $h=NF_L(p,F)$
\vskip 0.0cm \noindent
\hh \mbegin
\vskip 0.0cm \noindent
\hhh $h:=p$
\vskip 0.0cm \noindent
\hhh \mwhile exist $f\in F$ and a term $u$ of $h$ such that
\vskip 0.0cm \noindent
\hhhh $lm(f)|_L(u/cf(h,u))$\bb \mdo
\vskip 0.0cm \noindent
\hhhh \mchoose the first such $f$
\vskip 0.0cm \noindent
\hhhh $h:=h-(u/lt(f))f$
\vskip 0.0cm \noindent
\hhh \mend
\vskip 0.0cm \noindent
\hh \mend

\vskip 0.3cm
\noindent
{\em Correctness} and {\em termination} of this algorithm can be proved,
apparently, as they do for the conventional normal form
algorithm~\cite{Buch85,BWK93}.
Since involutive reductions form a fixed
subset of the conventional ones, generally, $NF_L(p,F)\neq NF(p,F)$.

\begin{definition}{\em A set $F$ is called {\em involutively autoreduced}
with respect to the given involutive division $L$, or {\em $L-$autoreduced},
if the set $lm(F)$ is $L-$autoreduced and every $f\in F$ has no terms
$t=cf(f,t)\,u\neq lt(f)$ with $cf(f,t)\neq 0$ and $u\in C_L(lm(F))$.
}
\end{definition}
Given an involutive division $L$ and a finite set $F$, the following
algorithm returns an $L-$autoreduced set
$H$, denoted by $H=Autoreduce_L(F)$, and such that $Id(F)=Id(H)$.

{\em Correctness} of the algorithm is obvious from the {\bf while-}loop
structure. Since the underlying set of involutive interreductions is a subset of the
conventional interreductions, its {\em termination} follows from that for
the conventional autoreduction~\cite{Buch85,BWK93}.

\vskip 0.3cm
\noindent
\hh Algorithm {\bf InvolutiveAutoreduction:}
\vskip 0.2cm
\noindent
\hh {\bf Input:}  $F$
\vskip 0.0cm \noindent
\hh {\bf Output:} $H=Autoreduce_L(F)$
\vskip 0.0cm \noindent
\hh \mbegin
\vskip 0.0cm \noindent
\hhh $H:=F$
\vskip 0.0cm \noindent
\hhh \mwhile exist $h\in H$ and $g\in H\setminus \{h\}$
\vskip 0.0cm \noindent
\hhhh such that $h$ is reducible modulo $g$\bb \mdo
\vskip 0.0cm \noindent
\hhhh \mchoose the first such $h$
\vskip 0.0cm \noindent
\hhhh $H':=H\setminus \{h\}$
\vskip 0.0cm \noindent
\hhhh $h':=NF_L(h,H)$
\vskip 0.0cm \noindent
\hhhh \mif $h'=0$\bb \mthen $H:=H'$
\vskip 0.0cm \noindent
\hhhh \melse $H:=H'\cup \{h'\}$
\vskip 0.0cm \noindent
\hhh \mend
\vskip 0.0cm \noindent
\hh \mend
\vskip 0.3cm
\noindent

\begin{theorem} If set $F$ is $L-$autoreduced,
then $NF_L(p,F)=0$ if and
only if $p$ is presented in terms of a finite sum of the form
\begin{equation}
p\in \S_F\subset \R\,,\quad \S_F=\{\ \sum_{ij} f_i\times u_{ij}\ |\ f_i\in F\,,
u_{ij}\in \T\ \} \label{nfzero}
\end{equation}
with $lm(u_{ij})\neq lm(u_{ik})$ for $j\neq k$.
\label{th_nfzero}
\end{theorem}

\vskip 0.2cm
\noindent
{\bf Proof}\ \
$\Longrightarrow :$ If $NF_L(p,F)=0$, then, by Definition~\ref{inv_red} of
involutive reductions, at each intermediate reduction step the current
value $p'$ of $p$ is rewritten as $p'\rightarrow p''=p'-f_i\times u_{ij}$.
Since the reduction chain is finite by admissibility of an ordering $\succ$,
the representation (\ref{nfzero}) holds.

$\Longleftarrow :$ Let $p$ is given by expression (\ref{nfzero}). Firstly,
we show that $lm(p)$ has an involutive divisor in the set $lm(F)$. For
this purpose select the leading term in the right hand side of (\ref{nfzero}).
It has
the form $s=lt(f_i\times u_{ij})=lt(f_i)\times u_{ij}$ with some $i,j$ and
cannot appear in any other term $lt(f_k)\times u_{kl}$. Otherwise,
the underlying monomial $s/lc(s)$ would have two involutive divisors
$lm(f_i)$ and $lm(f_k)$ what, by Proposition~\ref{unique_id},
would contradict the involutive autoreduction of $F$. Secondly, since $p$ is
involutively reducible, after each reduction step the representation
(\ref{nfzero}), obviously, still holds
providing the further reductions until the chain stops when we obtain
zero at a certain step. It just means that $NF_L(p,F)=0$.
\hfill{\Box}

\begin{corollary} If set $F$ is
$L-$autoreduced, then the $L-$normal form, for an arbitrary
algorithm of its computation and for
any polynomials $p_1,p_2$ and $p$, has the properties:{\em
\begin{enumerate}
\renewcommand{\theenumi}{(\roman{enumi})}
\item {\em Uniqueness: if $h_1=NF_L(p,F)$ and $h_2=NF_L(p,F)$ then
$h_1=h_2$.}
\item {\em Linearity: $NF_L(p_1+p_2,F)=NF_L(p_1,F)+NF_L(p_2,F)\,.$}
\end{enumerate} \label{pr_inf}
}
\end{corollary}

\vskip 0.2cm
\noindent
{\bf Proof}\ \
(i) By an involutive normal form algorithm,
$h_1=p-\sum_{ij}{f_i}\times u_{ij}$ and $h_2=p-\sum_{ij}{f_i}\times v_{ij}$.
Therefore, $h_1-h_2$ has the
representation (\ref{nfzero}), and $NF_L(h_1-h_2,F)=0$ by
Theorem~\ref{th_nfzero}. On the other hand, since $h_1$ and $h_2$ are
normal forms, they have no involutive divisors and so does $h_1-h_2$.
Hence, we have $h_1=h_2$.

(ii) Denote $p_1+p_2$ by $p_3$ and let
$$h_1=NF_L(p_1,F)\,,\quad h_2=NF_L(p_2,F)\,,\quad h_3=NF_L(p_3,F)\,.$$
Then $NF_L(h_3-h_1-h_2,F)=h_3-h_1-h_2$, since none of
$h_1,h_2,h_3$ has involutive divisors in $lm(F)$. In addition, because
$h_k=p_k-\sum_{ij}{f_i}\times v_{k;ij}$ ($k=1,2,3$), we have $h_3-h_1-h_2\in \S_F$. Thus,
by Theorem~\ref{th_nfzero}, $NF_L(h_3-h_1-h_2,F)=0$, and, hence, $h_3=h_1+h_2$.
\hfill{\Box}

\section{Involutivity Conditions}

\begin{definition}{\em  Multiplication of a polynomial $f\in F$ by a
 variable $x$ is called the {\em prolongation} of $f$. Given an involutive
 division specified by the set $lm(F)$, the prolongation
 is called {\em multiplicative} if $x$ is multiplicative for $lm(f)$, and
 {\em non-multiplicative}, otherwise.}
\end{definition}

\begin{definition}{\em An $L-$autoreduced set $F$ is called
{\em ($L-$)involutive} basis of $Id(F)$ if
\begin{equation}
(\forall f\in F)\ (\forall u\in \M)\ \ [\ NF_L(fu,F)=0\ ]\,. \label{iset}
\end{equation} \label{def_ibasis}
}
\end{definition}

\begin{proposition} Let $F$ be an involutive polynomial basis. Then
the monomial set $lm(F)$ is also involutive. \label{pr_ims}
\end{proposition}

\vskip 0.2cm
\noindent
{\bf Proof}\ \  It follows immediately from Definitions~\ref{inv_mset},
and~\ref{def_ibasis}
\hfill{\Box}

\vskip 0.2cm
\noindent
It is clear from Definition~\ref{def_ibasis} and the linearity of the
involutive
normal form, by Corollary~\ref{pr_inf}, that an {\em involutive basis}
provides
decision of the ideal membership problem. Hence, we have the following
corollary.

\begin{corollary} If set $F$ is $L-$involutive, then $p\in Id(F)$
if and only if $NF_L(p,F)=0$. In this case, obviously, the equality
$\S_F=Id(F)$
holds. \label{cr_ibasis}
\end{corollary}
The definition of involutive polynomial sets is the direct extension
of that for involutive monomial sets in Sect.4. The theorem below imparts
the constructive characterization of involutivity, which is the heart of the
involutive algorithms.

\begin{theorem}
 An $L-$autoreduced set $F$ is involutive
 with respect to a continuous involutive
 division $L$ if and only if the following conditions of local involutivity
 hold
\begin{equation}
 (\forall f\in F)\ (\forall x_i\in NM_L(lm(f),lm(F)))\
 \ [\ NF_L(f\cdot x_i,F)=0\ ]\,. \label{inv_cond}
\end{equation}
\label{th_inv_cond}
\end{theorem}

\vskip 0.2cm
\noindent
{\bf Proof}\ \
$\Longrightarrow :$ Since $x_i\in \M$ we are done.

$\Longleftarrow :$
An immediate consequence of~(\ref{inv_cond}) is local involutivity
of the set $lm(F)$ in accordance with Definition~\ref{mon_loc_inv}.
Then, by continuity of division $L$, this set is involutive. Thus,
for any $f\in F$ and any $u\in \M$ the monomial $lm(f)\cdot u$ has the
involutive divisor $lm(g)$, $g\in F$.

We claim that the polynomial $f\cdot u$ can be presented as
follows
\begin{equation}
f\cdot u=g\times v + \sum_{ij}f_i v_{ij}\,, \label{icf}
\end{equation}
where $v,v_{ij}\in \T$, $f_i\in F$ and relation
$lm(f\cdot u)=lm(g\times v) \succ lm(f_iv_{ij})$
holds for any term of the sum. Indeed, if $u$ is
multiplicative for $f$ we are trivially done. Otherwise $u$ contains
$x_k\in NM_L(f,lm(F))$. Then, the
local involutivity of $F$, by Theorem~\ref{th_nfzero}, yields the
representation
\begin{equation}
f\cdot x_k=g_1\times u_1 + \sum_{ij}f_i \times u_{ij} \label{li_rep}
\end{equation}
with $g_1\in F$ and $lm(f\cdot x_k)=lm(g_1u_1)\succ f_iu_{ij}$ for any
term under the summation sign. If monomial $u/x_k$ is
multiplicative for $g_1$, then~(\ref{icf}) immediately follows
from~(\ref{li_rep}) with $g=g_1$
and $v=u_1u/x_k$. Otherwise, multiply both sides of~(\ref{li_rep}) by $u/x_k$,
take a variable $x_m\in NM_L(g_1,lm(F))$, which is contained in $u/x_k$, and
apply the local involutivity conditions for $g_1\cdot x_m$. It gives
the relation
\begin{equation}
f\cdot u=(g_2\times u_2)u_1u/(x_kx_m) + \sum_{ij}f_i \tilde{u}_{ij}
\label{rep_uv}
\end{equation}
where inequality $lm(g_2)uu_1u_2/(x_kx_m)\succ lm(f_i\tilde{u}_{ij})$
holds for all $i,j$. If $uu_1/(x_kx_m)$ is still
non-multiplicative for $g_2$ the relation~(\ref{rep_uv}) can be further
rewritten by using the local involutivity conditions until we obtain
relation~(\ref{icf}). This is guaranteed by continuity of involutive
division $L$, because all the polynomials $g_1,g_2,\ldots \in F$ are
distinct, since their leading monomials, by construction, form the sequence
satisfying~(\ref{cont_cond}).

Next, similar rewriting the every
term $f_iv_{ij}$ in~(\ref{icf}) gives
$
f_i v_{ij}=f_k\times w_k + \sum_{lm}f_l w_{lm}
$
with
$lm(f_iv_{ij})=lm(f_k\times w_k) \succ lm(f_lw_{lm})$.
Proceeding with this way, by admissibility of ordering $\prec$, we
find, in a finite number of steps, that $f\cdot u\in \S_F$.
\hfill{\Box}

\vskip 0.2cm
\noindent
The next definition of partial involutivity is useful for the algorithmic
construction of involutive bases as we show below.

\begin{definition}{\em Given $v\in \M$ and an $L-$autoreduced set $F$, if
there exist $f\in F$ such that $lm(f)\prec v$ and
\begin{equation}
(\forall f\in F)\ (\forall u\in \M)\ (lm(f)\cdot u\prec v)\ \
[\ NF_L(fu,F)=0\ ]\,, \label{piset}
\end{equation}
then $F$ is called {\em partially involutive up to the monomial $v$}
with respect to the admissible ordering $\prec$. $F$ is
still said to be partially involutive up to $v$ if $v\prec lm(f)$ for
all $f\in F$.
}
\end{definition}
Looking at the proofs of Theorems~\ref{th_cont} and \ref{th_inv_cond} it is
easy to see that they prove also the following {\em conditions of partial
involutivity}.

\begin{corollary} Given a continuous involutive division $L$, an
 $L-$autoreduced set $F$ is partially involutive up to the monomial $v$ if
 and only if
\begin{equation}
(\forall f\in F)\ (\forall x_i\in NM_L(lm(f),lm(F)))\ (lm(f)\cdot x_i\prec v)\ \
[\ NF_L(f\cdot x_i,F)=0\ ]\,. \label{pinv_cond}
\end{equation}
\end{corollary}

\section{\Gr Bases and Involutive Bases}
\noindent
In~\cite{ZB93,ZB94} it was shown that a {\em Pommaret basis},
that is, involutive basis for Pommaret division, is also a \Gr basis,
though, generally,
not the reduced one. A similar property of a Janet basis was noticed
in~\cite{Zharkov94b}. The following theorem shows that such a relation
holds for any involutive division.

\begin{theorem} If set $F$ is $L-$involutive, then the equality of
the conventional and $L-$normal forms
\begin{equation}
(\forall p\in \R)\ \ [\ NF(p,F)=NF_L(p,F)\ ] \label{nf}
\end{equation}
holds for any normal form algorithm. \label{th_nf}
\end{theorem}

\vskip 0.2cm
\noindent
{\bf Proof}\ \  To prove the theorem it is sufficient to show
that any polynomial $p$ is reducible modulo $F$ if and only if it is
involutively reducible. But the latter statement is an easy consequence
of Definitions~\ref{inv_div} or \ref{3.1a} and~\ref{def_ibasis}. Indeed, if
$p$ is involutively reducible, then it is conventionally reducible.
Conversely, let the term $u$ have a divisor among the leading monomials of $F$,
that is, $u=lc(u)\,lm(f)\cdot v$ for some $f\in F$ and $v\in \M$. By
the condition (\ref{iset}) and Theorem~\ref{th_nfzero}, it implies
$f\cdot v=\sum_{ij}f_i\times u_{ij}$. Hence, $u$ has also the involutive
divisor in $lm(F)$. It is just that $f_i$ which satisfies the condition
$lm(f_i)\times u_{ij}=lm(f)\cdot v$ and is unique.
\hfill{\Box}

\begin{corollary} An involutive basis is a \Gr basis. \label{cr_igb}
\end{corollary}

\vskip 0.2cm
\noindent
{\bf Proof}\ \  According to the algorithmic characterization of \Gr
 bases~\cite{Buch65,Buch85,BWK93}
consider the S-polynomial of $f_i,f_j\in F$
\begin{equation}
S(f_i,f_j)=\frac{lcm(f_i,f_j)}{lt(f_i)}f_i - \frac{lcm(f_i,f_j)}{lt(f_j)}f_j\,.
\label{spol}
\end{equation}
From $S(f_i,f_j)\in Id(F)$, Corollary~\ref{cr_ibasis} and
Theorem~\ref{th_nf}, we have $ NF(S(f_i,f_j),F)=0$.

\hfill{\Box}

\begin{corollary} If set $F$ is partially involutive up to the monomial
$v$, then
\begin{equation}
(\forall p\in \R)\ (lm(p)\prec v)\ \ [\ NF(p,F)=NF_L(p,F)\ ]\,. \label{pnf}
\end{equation} \label{cr_pnf}
\end{corollary}

\vskip 0.2cm
\noindent
{\bf Proof}\ \  It follows by perfect analogy to the proof of
Theorem~\ref{th_nf}.
\hfill{\Box}

\vskip 0.2cm
\noindent
Note that while a Pommaret basis, if it exists for the given ideal,
is unique~\cite{ZB94}, this may not hold for other
involutive divisions. We demonstrate it by the following explicit
example.
\begin{example}{\em Two lexicographical ($x\succ y$) Janet bases $F_1$ and
$F_2$
\begin{eqnarray*}
F_1&=&\{xy^3-y,\overbrace{xy^2-1,xy-y^2,x-y}^{y},\overbrace{y^3-1}^x\}\,, \\[0.1cm]
F_2&=&\{x^2y^3-y^2,\overbrace{x^2y^2-y,x^2y-1,x^2-y^2}^y,
\overbrace{xy^3-y}^x,\overbrace{xy^2-1,xy-y^2,x-y}^{x,y},
\overbrace{y^3-1}^x\}\,,
\end{eqnarray*}
\vskip 0.1cm
\noindent
with indicated non-multiplicative variables, are involutive. It can easily
be verified. Both of them generate, obviously, the same ideal with the
\Gr basis $(x-y,y^3-1)$, which is also a Janet basis and, in this particular
case, coincides with the Pommaret basis.
} \label{exm_4}
\end{example}
As it was shown in Sect.4, given a polynomial set $F$ and an arbitrary
involutive division, the ideal $Id(F)$ may not have a finite
involutive basis.
For example, while a finite Pommaret basis exists for any
zero-dimensional ideal~\cite{Pommaret78,ZB94,Apel},
it may not exist
for a positive dimensional one. Generally, for positive dimensional
ideals, the existence of finite Pommaret basis can be achieved by
means of an appropriate linear transformation of
variables~\cite{Pommaret78,Apel}.

On the other hand, a noetherian involutive division, for
example, a Thomas or Janet one, implies the existence of finite
involutive bases for any polynomial ideals as the following
proposition shows.

\begin{proposition} If involutive division $L$ is noetherian, then any
polynomial ideal $Id(F)$ has a finite $L-$involutive basis.
\end{proposition}

\vskip 0.2cm
\noindent
{\bf Proof}\ \ Let $G$ be the reduced \Gr basis of $Id(F)$
which is finite for any polynomial ideal~\cite{Buch85,BWK93}.
If set $G$ is not involutive, then
complete it by non-multiplicative prolongations of its elements
just as it done in algorithm {\bf InvolutiveCompletion}. This means
that at every step of the completion we select a non-multiplicative
prolongation with the lowest leading term which is $L-$irreducible
modulo the current leading monomial set. By noetherity of $L$, in a
finite number of steps, a polynomial set $\tilde{G}$ will be
produced such that $lm(\tilde{G})$ be an $L-$autoreduced involutive
completion of $lm(G)$. Finally, $L-$autoreduction of the tales in
$\tilde{G}$ will give an $L$-involutive basis of $Id(F)$.
\hfill{\Box}

\section{Basic Algorithm}
\noindent
In this section we describe an algorithm for the construction of an
involutive basis. The algorithm is an improved version of one
presented in~\cite{ZB94} for Pommaret division and generalized to
any continuous noetherian division $L$ and any admissible ordering
$\succ$. The main optimization is based on the use of Buchberger's
chain criterion for avoiding unnecessary reductions introduced
in~\cite{Buch79} (see also~\cite{Buch85,BWK93}).

Corollary~\ref{cr_pnf} shows that for
any S-polynomial $S(f_i,f_j)$, given by formula (\ref{spol}),
both its conventional and
$L-$normal forms are vanishing as soon as the
conditions (\ref{pinv_cond}) are satisfied
up to the monomial $lcm(f_i,f_j)$.
According to Theorem~\ref{th_nfzero} and Corollary~\ref{pr_inf} the
conditions (\ref{pinv_cond}) can be
presented as $NF_L(S_L(f_i,f_j),F)=0\,,$
where $S_L(f_i,f_j)$
are just ($L-$involutive) S-polynomials of the special form
\begin{equation}
  S_L(f_i,f_j)=f_i\cdot x - f_j\times u_{jk}\,. \label{ispol}
\end{equation}
The following theorem gives the involutive form of Buchberger's chain
criterion.

\begin{theorem}
Let $F$ be a finite $L-$autoreduced polynomial set, and let
$g\cdot x$ be a non-multiplicative prolongation of $g\in F$.
Then $NF_L(g\cdot x,F)=0$ if the following holds
\begin{equation}
(\forall h\in F)\ (\forall u\in \M)\ \bigl(\ lm(h)\cdot u\prec
lm(g\cdot x)\  \bigr)\ \
[\ NF_L(h\cdot u,F)=0\ ]\,, \label{part_inv_cond}
\end{equation}
\begin{equation}
     (\exists f,f_0,g_0\in F)
\left[
\begin{array}{l}
 lm(f_0)|lm(f)\,,\ lm(g_0)|lm(g) \\[0.1cm]
 lm(f)|_L lm(g\cdot x)\,,\ lcm(f_0,g_0)\prec lm(g\cdot x) \\[0.1cm]
 NF_L\bigl(f_0\cdot \frac{lt(f)}{lt(f_0)},F\bigl)=
 NF_L\bigl(g_0\cdot \frac{lt(g)}{lt(g_0)},F\bigl)=0
\end{array}
\right]\,.
\label{inv_crit}
\end{equation}
\label{th_criteria}
\end{theorem}

\vskip 0.2cm
\noindent
{\bf Proof}\ \
Condition~(\ref{inv_crit}) yields that at least
one of polynomials $f,g$ can be considered as derived from $f_0,g_0$
by prolongations with at least one non-multiplicative among them.
If, for example, $lm(f_0)\neq lm(f)$, it leads to the equality
$f=f_0\cdot (lm(f)/lm(f_0))$ modulo $F$.

Thus, if the condition~(\ref{inv_crit}) holds, there
is a chain of polynomials in $F$ of the form
\begin{equation}
 f\equiv f_k,f_{k-1},\ldots,f_0,g_0,\ldots,g_{m-1},g_m\equiv g\,,
 \label{sp_chain}
\end{equation}
where $k+m>0$. Here $f$ or $g$ or both of them are produced by
prolongations, including non-multiplicative ones, of the polynomials
$f_i$ or $g_j$ in the chain whose indices are less than $k$ or $m$,
respectively.

The chain~(\ref{sp_chain}) has the property
$$ NF(S_L(f,f_{k-1}),F)=\cdots =NF(S(f_0,g_0),F)=\cdots =
NF(S_L(g_{m-1},g),F)=0\,.$$
This property is resulted from the observations as follow.
Consider relation
\begin{equation}
 lm(g)\cdot x=lm(f)\times w, \label{lt}
\end{equation}
which means that $w$ does not contain $x$. Otherwise, $g$ would be
reducible by $f$, and, hence, $F$ could not be $L-$autoreduced.
Thus, $lcm(f,g)=lm(g)\cdot x$. By admissibility of the
monomial ordering $\prec$, the least common multiple of the
leading monomials for pair of the neighboring polynomials in
the chain~(\ref{sp_chain}) is less than or equal to $g\cdot x$.
Then the above property of the chain follows immediately
from partial involutivity~(\ref{part_inv_cond}) of $F$ and
Corollary~\ref{cr_pnf}. Furthermore,
conditions~(\ref{part_inv_cond}-\ref{inv_crit}) imply
$NF_L(S(f_0,g_0),F)=NF(S(f_0,g_0),F)=0$, and
$NF_L(S_L(f_i,f_{i-1}),F)=NF(S(f_i,f_{i-1}),F)=0$ $(1\leq i\leq k)$
as well as $NF_L(S_L(g_{i-1},g_i),F)=NF(S(g_{i-1},g_i),F)=0$
$(1\leq i\leq m)$.

By construction,
$lcm(f,\ldots,f_1,f_0,g_0,g_1,\ldots,g)=lcm(f,g)$ what
leads~\cite{BWK93} to the representation
$S(f,g)= \sum_{ij} f_i u_{ij}$ where $f_i\in F$ and $lm(f_iu_{ij})\prec
lcm(f,g)=lm(g)\cdot x$. Then, condition~(\ref{part_inv_cond}), by
Corollaries~\ref{pr_inf} and~\ref{cr_pnf}, yields
$$NF_L(S_L(f,g),F)=NF(S(f,g),F)=0$$ in accordance with~\cite{Buch85,Buch79}.
\hfill{\Box}

\vskip 0.2cm
Before analysis of correctness and termination of the below algorithm, we give
some necessary clarifications.

First of all, the conventional autoreduction of the initial polynomial
set is done. It removes, in particular, all the predecessors of every
polynomial from the initial set.

Set $T$ collects all the triples $(g,u,P)$; $g$ is
an element in the current basis $G$; $u=lm(f)$ where $f\in G$ is the
predecessor of g, by a non-multiplicative prolongation of which
$g$ was derived, or $u=lm(g)$ if $g$ has no such predecessor in $G$;
$P$ is a set containing the non-multiplicative variables of $g$
have been used for its prolongations.

The current non-multiplicative prolongation $g\cdot x$ is selected to
be the lowest with respect to the ordering $\succ$. If there are several
different non-multiplicative prolongations with the same leading term, then
any of them may be selected. This selection strategy will be called
{\em normal}.

If the leading monomial of the current prolongation $g\cdot x$
is involutively reducible by the basis element $f\in G$, then
the other conditions in~(\ref{inv_crit}) are verified. The verification
is done in the form of comparison of $lcm(u,v)$ with $lcm(f,g)$,
where $u$ and $v$ are the second elements of the triples
containing $g$ and $f$, respectively. By Theorem~\ref{th_criteria},
the criterion~(\ref{inv_crit}) is false if and only if
$lcm(u,v)=lcm(f,g)=g\cdot x$. One should be also noted that
Buchberger's second criterion~\cite{Buch85} can be applied in the
involutive approach only in exceptional cases. Relation~(\ref{lt}) shows
that $lcm(f,g)=lm(f)lm(g)$ if and only if $lm(f)=x$ and $lm(g)=w$.

If the current prolongation is not reducible to zero, that is,
$h=NF_L(g\cdot x,G)\neq 0\,,$
then $h$ is added to $G$.

 After involutive autoreduction of the enlarged set $G$ an adjustment of
 the set $T$ is done. For an element $g\in G$ whose leading monomials was
 not mutually reduced, the second element $u$ in the triple is kept,
 if the leading term of the corresponding predecessor of $g$ was also
 not reduced. Otherwise, $u$ is replaced by its involutive divisor in
 $lm(G)$. Essentially new leading monomials, that is, those not multiple
 of any others occurring in $T$ before the autoreduction, are included in
 the refreshed $T$ with their actual leading monomials as the second
 elements of the triples.

\vskip 0.3cm
\noindent
\hh Algorithm {\bf InvolutiveBasis:}
\vskip 0.2cm
\noindent
\hh {\bf Input:}  $F$, a finite polynomial set
\vskip 0.0cm \noindent
\hh {\bf Output:} $G$, an involutive basis of the ideal $Id(F)$
\vskip 0.0cm \noindent
\hh \mbegin
\vskip 0.0cm \noindent
\hhh $G:=Autoreduce(F)$
\vskip 0.0cm \noindent
\hhh $T:=\emptyset$
\vskip 0.0cm \noindent
\hhh \mfore $g\in G$\bb \mdo $T:=T\cup \{(g,lm(g),\emptyset)\}$
\vskip 0.0cm \noindent
\hhh \mwhile exist $(g,u,P)\in T$ and $x\in NM_L(lm(g),lm(G))$
   \bb \mdo
\vskip 0.0cm \noindent
\hhhh \mchoose such $(g,u,P)$ and $x$ with the lowest $lm(g)\cdot x$
\vskip 0.0cm \noindent
\hhhh $T:=T\setminus \{(g,u,P)\} \cup \{(g,u,P\cup \{x\})\}$
\vskip 0.0cm \noindent
\hhhh \mif exist $f\in (f,v,D)\in T$ such that
 $lm(f)\,|_L\, lm(g\cdot x)$\bb \mthen
\vskip 0.0cm \noindent
\hhhhh \mif $lcm(u,v) = lm(g)\cdot x$\bb \mthen $h:=NF_L(g\cdot x,G)$
\vskip 0.0cm \noindent
\hhhhhh \mif $h\neq 0$\bb \mthen $T:=T\cup \{(h,lm(h),\emptyset)\}$
\vskip 0.0cm \noindent
\hhhh \melse
 $h:=NF_L(g\cdot x,G)$
\vskip 0.0cm \noindent
\hhhh $T:=T\cup \{(h,u,\emptyset)\}$
\vskip 0.0cm \noindent
\hhhh $G:=Autoreduce_L(G\cup \{h\})$
\vskip 0.0cm \noindent
\hhhh $Q:=T$
\vskip 0.0cm \noindent
\hhhh $T:=\emptyset$
\vskip 0.0cm \noindent
\hhhh \mfore $g\in G$\bb \mdo
\vskip 0.0cm \noindent
\hhhhh \mif exist $(f,u,P)\in Q$ such that $lm(f)=lm(g)$\bb \mthen
\vskip 0.0cm \noindent
\hhhhhh \mchoose $g_1\in G$ such that $lm(g_1)|_Lu$
\vskip 0.0cm \noindent
\hhhhhh $T:=T\cup \{(g,lm(g_1),P)\}$
\vskip 0.0cm \noindent
\hhhhh \melse $T:=T\cup \{(g,lm(g),\emptyset)\}$
\vskip 0.0cm \noindent
\hhh \mend
\vskip 0.0cm \noindent
\hh \mend
\vskip 0.3cm

{\em Correctness}. As we have shown, criterion~(\ref{inv_crit}) is used
in algorithm {\bf InvolutiveBasis} in accordance with
Theorem~\ref{th_criteria}. It is easy to show that there
is the unique polynomial $g_1\in G$ which is chosen in the
inner {\bf for each}-loop such that $lm(g_1)$ involutively
divides $u$. Indeed, if the leading term of the predecessor $h$ of
$g$ with $u=lm(h)$ has not been reduced, then $g_1=h$. Otherwise,
there is $g_1\in G$ such that $g_1\neq h$ and $lm(g_1)|_Lu$.
The uniqueness of $g_1$ for the autoreduced set $G$ is an immediate
consequence of the property (v) in Definition~\ref{inv_div}.
Besides, the replacement of $u$ by $g_1$ does not violate, obviously,
the conditions for applicability of the criterion.
Furthermore, from Corollary~\ref{cr_pnf} it follows: (i) a leading
monomial, being involutively reducible at some step of the algorithm,
will never appear again among the leading monomials; (ii) there
is no need in recomputing zero reductions after enlargement of an
intermediate polynomial set.
This enables one to assign the set $P$ of the used non-multiplicative
variables for
polynomial $f$ to the corresponding polynomial $g$ with $lm(g)=lm(f)$
as it is done in the inner {\bf for each}-loop. Such an optimization
allows one to avoid the repeated prolongations.

Therefore, if division $L$ is continuous, and the algorithm terminates,
then it produces, by Theorem~\ref{th_inv_cond}, the involutive basis.
The termination holds if and only if the set $P$ in each triple
$(g,u,P)\in T$ contains all non-multiplicative variables for basis
element $g$. It just means that any non-multiplicative prolongation
of every element in $G$ is reduced to zero, and, hence, $G$ is involutive.

\vskip 0.2cm
{\em Termination.}  Note that the initial value of the leading
monomial set
$$U_0=lm(Autoreduce(F))$$ is determined by the input set
$F$ subjected to the conventional autoreduction. Since
only those monomials occur in the leading monomial set which
have not been reducible at some step of the algorithm, the change
in set $U=lm(G)$ after running the
{\bf while}-loop may take place only in two cases:
\begin{enumerate}
\renewcommand{\theenumi}{(\roman{enumi})}
\item $lm(g)\cdot x$ has no involutive divisors in $U$. In this case
 $U$ is enlarged to include $lm(g)\cdot x$.
\item $g\cdot x$ is reducible by elements of $U$. Then $U$ is enlarged
 to include $lm(h)$, where $h=NF_L(g\cdot x,G)\neq 0$ and $lm(h)$
 is not multiple,
 in the conventional sense, of any elements in $U_0$.
\end{enumerate}
 The number of different $lm(h)$ occurring in case (ii) is finite by
 Dickson's lemma (Becker, Weispfenning and Kredel, 1993).
 Recall also that  algorithms {\bf InvolutiveAutoreduction}
 and {\bf InvolutiveNormalForm} always terminate (Sect.5).

 Thus,
 the algorithm termination is determined by that of algorithm
 {\bf InvolutiveCompletion} considered in Sect.4. It follows that
 algorithm {\bf InvolutiveBasis} terminates for any noetherian division
 and arbitrary input polynomial set $F$. If $L$ is not noetherian,
 then termination may not hold if
 an intermediate set $U=lm(G)$ is not finitely generated with
 respect to $L$ as the below Example~\ref{exm_5} shows.
 In the case of Pommaret division the algorithm terminates, however,
 for any degree compatible
 ordering and any
 zero-dimensional ideal~\cite{ZB94}.
 Because the involutive division $L$ is continuous, once
 algorithm {\bf InvolutiveCompletion} terminates, an $L-$completion
 $\tilde{U}$ of $U$ will be constructed such that autoreduction of the
 corresponding set $G$ does not produce new leading monomials. $G$ is,
 obviously, the output involutive basis.

Proposition~\ref{pr_finite} implies, in particular, the
algorithm termination for Thomas and Janet divisions. However, for Pommaret
division, which is not noetherian, the algorithm may not terminate even in
the case when there is a finite Pommaret basis but the ordering is not
degree compatible as the following simple example shows.

\begin{example}{\em The set $F=\{x^2-1,xy-1,z\}$ generates a
zero-dimensional ideal with the lexicographical Pommaret basis
$(x\succ y\succ z)$ given by
$G=\{x-y,y^2-1,yz,z\}.$ However, following the above algorithm we have to
choose $z\cdot y$ as the first prolongation which is lexicographically
lowest. Since polynomial $h=yz$ has no Pommaret divisors among $lm(F)$,
we find $F\cup \{yz\}$ as an intermediate basis. The next lowest
prolongation is $yz\cdot y$ again has no Pommaret divisors among the
leading monomials of the enlarged set. Exploring this procedure further
produces the infinite involutively irreducible set
$$ \{x^2-1,xy-1,z,yz,y^2z,\ldots,y^kz,\ldots\}\quad k\in \N\,.$$
} \label{exm_5}
\end{example}
It is well-known~\cite{Pommaret78,ZB93,ZB94,Apel}
that positive dimensional ideals may not have finite Pommaret
bases. Example~\ref{exm_2} illustrates this fact at the monomial
level.
The following more non-trivial example shows the output of algorithm
{\bf InvolutiveBasis} for Pommaret and Janet divisions in the
case of polynomial ideal.
\begin{example}{\em Cyclic 4-th roots.
\vskip 0.3cm
\begin{center}
\bg {|r|r|l|} \hline\hline
$NM_J$          & $NM_P$ &\hspace*{1.0cm} Initial Polynomial Set  \\  \hline
$x_2$           & $-$      & $x_1+x_2+x_3+x_4$ \\
$x_3$           & $x_1$  & $x_1x_2+x_2x_3+x_3x_4+x_4x_1$    \\
$x_4$           &$ x_1,x_2$ & $x_1x_2x_3+x_2x_3x_4+x_3x_4x_1+x_4x_1x_2$   \\
 $-$ & $x_1,x_2,x_3$  &  $x_1x_2x_3x_4-1$  \\ \hline \hline
\nd
\end{center}
} \label{exm_6}
\end{example}
\vskip 0.3cm
Here we choose the degree-reverse-lexicographical-ordering with the order of
variables as in Sect.2. Note that, since the initial set is not autoreduced,
the inclusion $NM_J\subseteq NM_P$ (see Proposition~\ref{id_P_J})
does not hold.

Application of algorithm {\bf InvolutiveBasis} gives the following form
of Janet and Pommaret bases
\vskip 0.3cm
\begin{center}
\bg {|r|r|l|} \hline\hline
$NM_J$          & $NM_P$ &\hspace*{0.7cm} Janet and Pommaret Bases  \\
\hline
$ - $  & $-$    & $x_1+x_2+x_3+x_4$ \\
$x_1 $ & $x_1$  & $x_2^2 + 2x_2x_4 + x_4^2$    \\
$x_1,x_2$  &$ x_1,x_2$ & $x_2x_3^2 + x_3^2x_4 - x_2x_4^2 - x_4^3$   \\
$x_1,x_2,x_3$ & $x_1,x_2,x_3$ &  $x_2x_3x_4^2 + x_3^2x_4^2 - x_2x_4^3 +
  x_3x_4^3 - x_4^4 - 1$  \\
$x_1,x_2,x_3$ & $x_1,x_2,x_3$ & $x_2x_4^4 + x_4^5 - x_2 - x_4$  \\
$x_1,x_2,x_3$ & $x_1,x_2,x_3$ &
$x_3^2x_4^4 + x_2x_3 - x_2x_4 + x_3x_4 - 2x_4^2$   \\
$x_1,x_2 $  & $x_1,x_2,x_3$ & $x_3^3x_4^2 + x_3^2x_4^3 - x_3 - x_4$ \\ \hline
&  & \\[-0.3cm]
    & $x_1,x_2,x_3$ & $x_3^4x_4^2 + x_2x_3 - x_3^2 - x_2x_4 + x_3x_4 - x_4^2$  \\
  &                 & ........................................................... \\
  &                 &  \\
\hline \hline
\nd
\end{center}
\vskip 0.3cm
The Janet basis consists of the upper seven polynomials and coincides
with the
\Gr basis, while the Pommaret basis is infinite and contains also
prolongations of the seventh polynomial with respect to its
non-multiplicative variable $x_3$. Note that the ideal is
one-dimensional, that is why it does not have a finite Pommaret
basis.

The algorithm {\bf InvolutiveBasis} has been implemented in Reduce 3.5 for
the degree-reverse-lex\-i\-co\-graph\-i\-cal-ordering and Pommaret
division refined in a certain way to provide the algorithm termination
for any polynomial ideal. This refinement is equivalent to the dynamical
incorporation of some noetherian involutive division in the computational
process. Its detailed description will be given elsewhere. In addition,
the current package called INVBASE is considerably faster than
previous version~\cite{ZB94}, in particular, since it
uses the criterion~(\ref{inv_crit}).

\noindent
Experimentally, we observed much smoother behavior of the algorithm
{\bf Involutive
Basis} with respect to Buchberger algorithm\footnote{More precisely,
with respect to its implementation in Reduce 3.5.} as the ordering changes.
Consider, for instance, the following example.
\begin{example} {\em Cyclic 6-th roots.
\vskip 0.2cm
\noindent
\hspace*{0.5cm}$ x_1 + x_2 + x_3 + x_4 + x_5 + x_6\,,$ \\
\hspace*{0.5cm}$ x_1x_2 + x_2x_3 + x_3x_4 + x_4x_5 + x_5x_6 + x_6x_1\,,$ \\
\hspace*{0.5cm}$ x_1x_2x_3 + x_2x_3x_4 + x_3x_4x_5
    + x_4x_5x_6 + x_5x_6x_1 + x_6x_1x_2\,,$ \\
\hspace*{0.5cm}$x_1x_2x_3x_4 + x_2x_3x_4x_5 + x_3x_4x_5x_6
    + x_4x_5x_6x_1 + x_5x_6x_1x_2 + x_6x_1x_2x_3\,, $ \\
\hspace*{0.5cm}$ x_1x_2x_3x_4x_5 + x_2x_3x_4x_5x_6 + x_3x_4x_5x_6x_1
    + x_4x_5x_6x_1x_2 + x_5x_6x_1x_2x_3 + x_6x_1x_2x_3x_4\,, $ \\
\hspace*{0.5cm}$  x_1x_2x_3x_4x_5x_6 - 1\,. $}
\label{exm_7}
\end{example}
The next table gives the timings of INVBASE on an 66 Mhz MS-DOS based
AT/486 computer for different degree-reverse-lex\-i\-co\-graph\-i\-cal-orderings.
\vskip 0.3cm
\begin{center}
\bg {|c|c|c|} \hline\hline
\hspace*{0.5cm}Ordering\hspace*{0.5cm}&\hspace*{0.5cm}Timing (sec.)\hspace*{0.5cm} \\ \hline
$x_1\succ x_2\succ x_3\succ x_4\succ x_5\succ x_6$ & 1040 \\
$x_1\succ x_2\succ x_4\succ x_6\succ x_3\succ x_5$ & 514     \\
$x_1\succ x_2\succ x_4\succ x_6\succ x_5\succ x_3$ & 437     \\
$x_1\succ x_2\succ x_6\succ x_3\succ x_4\succ x_5$ & 1066     \\
$x_1\succ x_3\succ x_4\succ x_5\succ x_2\succ x_6$ & 604     \\
$x_1\succ x_3\succ x_4\succ x_6\succ x_5\succ x_2$ & 136     \\
$x_1\succ x_4\succ x_2\succ x_3\succ x_5\succ x_6$ & 993      \\
$x_1\succ x_4\succ x_5\succ x_6\succ x_2\succ x_3$ & 1001     \\
$x_1\succ x_5\succ x_3\succ x_4\succ x_6\succ x_2$ & 364     \\
$x_1\succ x_5\succ x_6\succ x_2\succ x_3\succ x_4$ & 1045     \\
$x_1\succ x_6\succ x_3\succ x_2\succ x_4\succ x_5$ & 1012     \\
$x_1\succ x_6\succ x_5\succ x_2\succ x_4\succ x_3$ & 590    \\
\hline \hline
\nd
\end{center}
\vskip 0.3cm
Comparison with the package GROEBNER implementing Buchberger algorithm
on the same Reduce 3.5 platform shows that its corresponding timings are not
only much larger than those presented in the table, but also vary
dramatically with the order of the variables. This fact was already observed
in~\cite{ZB94} where some comparative data for GROEBNER and
the previous version of the INVBASE package are presented.

\section{Conclusion}
\noindent
Buchberger algorithm and the involutive one are based on different
rewriting techniques, namely, on the use of S-polynomials and prolongations,
respectively, as well as on distinct reduction processes.
Nevertheless, as we demonstrate in
this paper, they are in fact very interconnected. If, as we
propose in the
algorithm {\bf InvolutiveBasis}, we choose the current prolongation in
increasing
order with respect to given monomial ordering, then the conventional and
involutive
normal form will coincide. What is more, the involutive
reduction of the prolongation is equivalent to the consideration of a
certain
S-polynomial. Just this fact makes it possible to use Buchberger's
criteria.

Recently another interesting facet of interrelation of both methods
was discovered by Apel~\cite{Apel}, namely, that Pommaret bases can
be associated with \Gr ones in appropriate graded structures. Earlier
such \Gr bases were intensively investigated in more general context by
Mora~\cite{Mora}. That
observation gives an opportunity to algorithmically construct
Pommaret bases
whenever they exist~\cite{Apel}.
Though such an analogy also enables one to take advantage of Buchberger's
criteria, it is restricted to Pommaret division.

Thus, all the above, as well as computer experiments with both techniques,
offers a clearer view of the most optimal computational procedures.

There is no question that any algorithmic improvement of the \Gr basis
and involutive techniques at the algebraic level has an analogous
optimization at the differential level, at least for linear partial
differential equations~\cite{Gerdt95}.

\section{Acknowledgements}
\noindent
The initial version of this paper was written on the eve of the
first anniversary of death of our mutual friend and coauthor
Alyosha Zharkov who made invaluable contribution to the formation
of the involutive approach to commutative algebra. We devote the
present work to his memory. The authors are grateful to J.Apel for
numerous fruitful remarks and comments. This work was supported in
part by the RFBR grant No. 96-15-96030.

\end{document}